\documentclass{amsart}
\usepackage{amsmath}
\usepackage{amssymb}
\usepackage{latexsym}
\usepackage{amsfonts}

\newtheorem{Pa}{Paper}[section]
\newtheorem{Tm}[Pa]{{\bf Theorem}}
\newtheorem{La}[Pa]{{\bf Lemma}}

\newcommand{\CC}{{\mathchoice
{\setbox0=\hbox{$\displaystyle\rm C$}\hbox{\hbox
to0pt{\kern0.4\wd0\vrule height0.9\ht0\hss}\box0}}
{\setbox0=\hbox{$\textstyle\rm C$}\hbox{\hbox
to0pt{\kern0.4\wd0\vrule height0.9\ht0\hss}\box0}}
{\setbox0=\hbox{$\scriptstyle\rm C$}\hbox{\hbox
to0pt{\kern0.4\wd0\vrule height0.9\ht0\hss}\box0}}
{\setbox0=\hbox{$\scriptscriptstyle\rm C$}\hbox{\hbox
to0pt{\kern0.4\wd0\vrule height0.9\ht0\hss}\box0}}}}

\begin{document}
\title[On reproducing kernel Hilbert spaces]
{On the reproducing kernel Hilbert spaces associated with the
fractional and bi--fractional Brownian motions}
\author[D. Alpay]{Daniel Alpay}
\address{Department of Mathematics \\
Ben-Gurion University of the Negev \\
P.O. Box 653 \\
84105 Beer-Sheva, Israel} \email{\tt dany@math.bgu.ac.il}
\author[D. Levanony]{David Levanony}
\address{Department of Electrical Engineering\\
Ben-Gurion University of the Negev \\
P.O. Box 653 \\
84105 Beer-Sheva, Israel} \email{\tt levanony@ee.bgu.ac.il}
\thanks{D. Alpay thanks the Earl Katz family for endowing the chair which
supports his research}

\begin{abstract}
We present decompositions of various positive kernels as
integrals or sums of positive kernels. Within this framework we
study the reproducing kernel Hilbert spaces associated with the
fractional and bi-fractional Brownian motions. As a tool, we
define a new function of two complex variables, which is a natural
generalization of the classical Gamma function for the setting we
consider.
\end{abstract}

\maketitle
%

\section{Introduction}
\setcounter{equation}{0}
In this paper we present expansions of the form
\begin{equation}
\label{La_Fourche}
{\mathbf K}(t,s)=\sum_{n=1}^\infty K^n(t,s)K_n(t,s)
\end{equation}
where $K(t,s)$ and $K_n(t,s),\, n=1,2,\ldots$ are positive kernels
on the real line, for a family of kernels ${\mathbf K}(t,s)$ of the
the form
\begin{equation}
\label{Maubert_Mutualite}
{\mathbf K}(t,s)=\varphi(r(t)+r(s))-\varphi(r(t-s))
\end{equation}
with appropriate choices of functions $r$ and $\varphi$. The
functions $r$ we consider are real--valued and such that
$r(0)=0$, but it is
better to assume them complex--valued and with 
a possibly non zero $r(0)$ at this stage of the
exposition. The case where $\varphi(x)=x$, corresponds to a class
of functions introduced by J. von Neumann and I. Schoenberg in
the late 1930s (see \cite{MR1503439}, \cite{MR0004644}), and M.G.
Krein  in the 1940s (see \cite{MR0012176} and the 1954 paper
\cite{krein_1954}; Krein's papers are reprinted in
\cite{MR1321817}). This class plays an important role both in analysis and
in the theory of stochastic processes. To recall the role of the
functions $r$ such that the kernel
\begin{equation}
\label{Duroc_Vaneau}
K_r(t,s)=r(t)+r(s)^*-r(t-s)-r(0),
\end{equation}
is positive on the real line, we first give a definition: A
{\em helical arc} is a continuous map $t\mapsto \xi_t$ from the real
line into a metric space (with metric $d$) such that
$d(\xi_t,\xi_s)$ is a function of $t-s$:
\begin{equation}
\label{160407} d(\xi_t,\xi_s)=
\rho(t-s).
\end{equation}
When the metric space is a {\sl real} Hilbert space,
it is a consequence
of a result of K. Menger (see \cite{menger}) that a function
$\rho$ satisfies \eqref{160407} if and only if the kernel 
\[
K_{\rho^2}(t,s)=\rho^2(t)+\rho^2(s)-\rho^2(t-s)
\]
is positive on the real line; see
\cite[Lemma 1, p. 229]{MR0004644}. Using this fact,
von Neumann and Schoenberg
give in \cite[Theorem 1, p. 229]{MR0004644}  an integral
representation of the functions $\rho$ corresponding to a helical arc
with values in a real Hilbert space. These are the
functions $\rho$ such that:
\begin{equation}
\label{eureka} \rho^2(t)=\int_0^\infty
\frac{\sin^2tu}{u^2}d\sigma(u),
\end{equation}
where the positive real measure $\sigma$ is such that
\[
\int_1^\infty\frac{d\sigma(u)}{u^2}<\infty.
\]
We note the formula (see \cite[(1.6) p. 229]{MR0004644})


\[
\rho^2(t)+\rho^2(s)-\rho^2(t-s)=\int_0^\infty\frac{
(1-\cos2tu)(1-\cos2su)+\sin2tu\sin 2su}{2u^2}d\sigma(u)
\]
and hence the kernel $\rho^2(t)+\rho^2(s)-\rho^2(t-s)$, which is
of the form  \eqref{Duroc_Vaneau}, is positive on the real line.
The result of von Neumann and Schoenberg expresses in particular
that the converse is true: if $r$ is an even real-valued function
such that $r(0)=0$ and $K_r(t,s)$ is positive on the real line, then
$r$ is
of the form \eqref{eureka}.\\

The motivation of von Neumann and Schoenberg
comes from imbedding problems of metric spaces into Hilbert
spaces; we will not address this question here, but refer the
interested reader to
\cite[pp.
81--84]{MR86b:43001} for further information.\\

In \cite{MR0012176}, Krein considers helical arcs defined on an
open interval $(a,b)\subset{\mathbb R}$, with values in a {\sl
complex} Hilbert space. Condition \eqref{160407} is replaced by
the condition on the inner product, that is, by the condition that the 
function
\[
B(t,s)=\langle \xi_{s+v}-\xi_v,\xi_{t+v}-\xi_v\rangle_{\mathcal
H}\]
does not depend on $v$. Krein states that a function $B(t,s)$
corresponds to a helical arc defined on the real line if and only
if it is of the form
\begin{equation}
B(t,s)=\int_{\mathbb R}
\frac{e^{itu}-1}{u}\frac{e^{-isu}-1}{u}dm(u)
\end{equation}
for an appropriate increasing function $m$; see \cite[Theorem
4]{MR0012176} and \cite[p. 115]{MR1321817}. One can then take
\[
\xi_t(u)=e^{itu}\quad{\rm and}\quad {\mathcal H}= {\mathbf
L}_2\left(\frac{dm(u)}{u^2}\right).
\]

In particular,
$B(t,s)$ defines a positive kernel on ${\mathbb R}$. Setting
\[
g(t)=i\gamma t+\int_{\mathbb
R}\left\{e^{itu}-1-\frac{itu}{u^2+1}\right\}\frac{dm(u)}{u^2}
\]
one has
\[
B(t,s)=g(t-s)-g(s)-g(s)^*.
\]
It follows (see also the paper \cite[equations (11.16) and
(11.17), p. 402]{KL_jot} of M.G. Krein and H. Langer and the
preprint \cite{denisov}) that a function $r(t)$, satisfying to
$r(-t)=r(t)^*$, is such that the kernel $K_r(t,s)$ in
\eqref{Duroc_Vaneau} is positive on the real line if and only if
it can be written as
\begin{equation}
\label{babylone}
r(t)=r_0+i\gamma t-\int_{\mathbb
R}\left\{e^{itu}-1-\frac{itu}{u^2+1}\right\}\frac{dm(u)}{u^2}
\end{equation}
where $r_0=r(0)$ and $\gamma$ are real numbers and where $m$ is a
positive measure on ${\mathbb R}$ such that
$$\int_{\mathbb
R}\frac{dm(u)}{u^2+1}<\infty.$$
Furthermore, equation \eqref{babylone} implies the representation
of the kernel
\begin{equation}
\label{marseille}
r(t)+r(s)^*-r(t-s)-r(0)=\int_{\mathbb
R}\frac{e^{itu}-1}{u}\frac{e^{-isu}-1}{u}dm(u)
\end{equation}
as an inner product.
See \cite[equation (11.17) p. 402]{KL_jot}.\\

Since $m$ is positive and hence real valued, we have:
\begin{equation}
\begin{split}
{\rm Re}\left\{r(t)\right\}&=\int_{\mathbb R} \dfrac{1-\cos(tu)}{u^2}dm(u)\\
{\rm Im}\left\{r(t)\right\}&=\int_{\mathbb R}
\left\{-\sin(tu)+\dfrac{tu}{u^2+1}\right\}dm(u)
\end{split}
\end{equation}
In particular, when the measure is even and $r(0)=0$, one has
\[
r(t)=\int_{\mathbb R} \dfrac{1-\cos(tu)}{u^2}dm(u).\]
With the equality $1-\cos 2u=2\sin^2 u$, and a change of variables,
one then recognizes formula \eqref{eureka} of von
Neumann and Schoenberg. Note also that the
kernel $K_r(t,s)$ takes the simpler form
\begin{equation}
\label{mathilde_est_revenue}
K_r(t,s)=r(t)+r(s)-r(t-s).
\end{equation}

The choice where $dm(u)=\frac{1}{\pi}|u|^{1-2H}du$ with
$H\in(0,1)$, corresponds to covariance function of the fractional
Brownian motion
with Hurst parameter $H$. Indeed, we have the following computation:
\[
\int_0^\infty\dfrac{1-\cos(tu)}{u^2}|u|^{1-2H}
du=\int_0^\infty\dfrac{1-\cos
y}{y^2t^{-2}} y^{1-2H} t^{-(1-2H)}\dfrac{dy}{t}\\
=c({1-2H})t^{2H},
\]
where
\[
c({1-2H})=\int_0^\infty (1-\cos
y)y^{-1-2H}dy=\dfrac{1}{(1-2H){2H}}\int_0^\infty y^{1-2H}\cos y
dy,
\]
where we have assumed that $H\not=1/2$, and where we have
integrated twice by parts. The last integral is computed by
integrating the function
\[f(z)=e^{iz}e^{(1-2H) \ln z},\]
where $\ln z$ is the principal value of the logarithm in ${\mathbb
C}\setminus{\mathbb R}_-$, along an appropriate contour and is
found equal to
\[
\Gamma(2-2H)\cos(H\pi),
\]
where $\Gamma$ denotes Euler's Gamma function.
Thus, \(r(t)=\dfrac{V_H}{2}|t|^{2H}\), where
\[
V_H=\dfrac{\Gamma(2-2H)\cos(H\pi)}{\pi (1-2H)H}.
\]

The function
\begin{equation}
\label{ioanna}
K_{H,\alpha}(t,s)=(|t|^{2H}+|s|^{2H})^\alpha-|t-s|^{2H\alpha},\
\end{equation}
with
\[
0<\alpha\le 1\quad{\rm and}\quad 0<H< 1,
\]
is a special case of \eqref{Maubert_Mutualite} with
$r(t)=|t|^{2H}$ and $\varphi(x)=x^\alpha$. It was introduced in
\cite[Exercise 2.12 (h) p. 79]{MR86b:43001}, is positive in
${\mathbb R}$, and, by Kolmogorov's theorem, may assume the role
of the covariance function of a zero mean Gaussian stochastic
process. This process was introduced and first studied by Houdr\'e
and Villa in \cite{MR2037165}, who coined for it the term {\sl
bi--fractional Brownian motion}. The classical Brownian motion
corresponds to the choice $\alpha=1$ and $H=1/2$ while the choice
$\alpha=1$ and $H<1$ ($H\not =1/2$) corresponds, up to a
multiplicative constant factor, to the covariance function of the
fractional Brownian motion. For $\alpha<1$, one can view this
kernel as a generalization of the covariance function of the
fractional Brownian motion. It has been recently the topic of a
number of
papers; see \cite{MR2037165}, \cite{russo_tudor}.\\

It is well known that the reproducing kernel space, associated
with the covariance function of the Brownian motion, is the
Sobolev space of functions absolutely continuous on the real line
and with square summable derivatives; see \cite[p. 25]{Malliavin},
\cite{MR99f:60082}. The reproducing kernel Hilbert space
associated with the covariance function of the fractional Brownian
motion has also been characterized in a number of places; see for
instance \cite[Theorem 4.1 p. 948]{fbm88}.\\

We note that Schoenberg \cite[(4), p. 788]{MR1503370} also considers the
multidimensional kernel case.\\

We now turn to the outline of the paper. The paper consists of
five sections besides the introduction. Sections 2 and 3 are of a
review nature. The main results appear in Sections 4,5 and 6. In
Section 2 we recall some of M.G. Krein's results on positive
kernels of the form \eqref{Duroc_Vaneau} and of the associated
reproducing kernel Hilbert spaces. These results have been
rediscovered within the setting of the fractional Brownian motion
a number of times without reference to Krein's work. In the third
section we recall several facts on integrals and products of
positive kernels and of the associated reproducing Hilbert kernel
spaces. Section 4 contains two proofs of the positivity of the
kernel of the fractional Brownian motion. The first is valid for
any $H\in(0,1)$ while the second one is valid only for
$H\in(0,1/2)$. 
In both instances, the positivity of
the kernel is proved by exhibiting it as a sum of positive
kernels of the form \eqref{La_Fourche}; see formulas
\eqref{total-recall-nice-1991} and \eqref{ohohoh} below. In
Section 5 we consider the case of the bi-fractional Brownian
motion and associate to its covariance function an expansion of
the form \eqref{La_Fourche}. In Section 6 we study some nonlinear
transforms associated with the expansions \eqref{La_Fourche}.\\

A note on terminology: Throughout we use the terms
positive kernels, positive functions, covariance kernels or
covariance functions interchangeably.

\section{The reproducing kernel Hilbert space ${\mathcal H}(K_r)$}
\setcounter{equation}{0}
Recall that a function $K(t,s)$ defined on a set $\Omega$ is said
to be positive if it is hermitian:
\[K(t,s)=K(s,t)^*,\quad \forall t,s\in \Omega,\]
and if, for every integer $n$ and every choice of points
$t_1,\ldots, t_n\in\Omega$, the $n\times n$ matrix with $\ell j$
entry $K(t_\ell,t_j)$ is nonnegative. In preparation for the
forthcoming sections (see in particular formula \eqref{St
Philippe du Roule, ligne 9} below), one also may recall that sums
and products of positive functions are still positive functions.
The Aronszajn--Moore theorem states that a function $K(t,s)$ is
positive on a set $\Omega$ if and only if it can be written as
\begin{equation}
\label{zxcvbnm}
K(t,s)=\langle f_s,f_t\rangle_{\mathcal C},
\end{equation}
where ${\mathcal C}$ is a auxiliary Hilbert space and where
$s\mapsto f_s$ is a function from $\Omega$ into ${\mathcal C}$.
The space ${\mathcal C}$ can be chosen to be the reproducing
kernel Hilbert space associated with $K(t,s)$. This is a Hilbert
space of functions on $\Omega$, denoted by ${\mathcal H}(K)$,
uniquely determined by the following two conditions: $(a)$ For
every $s\in\Omega$ the function
\[K_s\,\,:\,\, u\mapsto K(u,s)\]
belongs to ${\mathcal H}(K)$, and: $(b)$ for every $s\in\Omega$
and $F\in{\mathcal H}(K)$
\begin{equation}
\langle F, K_s\rangle_{{\mathcal H}(K)}=F(s). \label{rk_pr}
\end{equation}

Setting $f_s=K_s$ in \eqref{rk_pr}, leads to the factorization
\eqref{zxcvbnm}.\\

If $K$is positive on $\Omega$ so is
the function $mK$, where $m$ is a strictly positive function. We
recall that
\begin{equation}
\label{poisse}
\|F\|^2_{{\mathcal H}(K)}=
m\|F\|^2_{{\mathcal H}(mK)}.
\end{equation}
For future reference we now recall the
characterization of the reproducing kernel space when a
representation \eqref{zxcvbnm} is available.
\begin{La}
\label{Vaugirard}
Let $K(t,s)$ be of the form \eqref{zxcvbnm} and let ${\mathcal
C}(f)$ be the closed linear span in ${\mathcal C}$ of the functions
$f_s$ where $s$ runs throughout $\Omega$. Then the reproducing kernel
Hilbert space with reproducing kernel $K(t,s)$ is the set of
functions of the form
\[
F(t)=\langle x, f_t\rangle_{\mathcal C},\quad x\in{\mathcal C}(f)
\]
with norm
\[
\|F\|_{{\mathcal H}(K)}=\|x\|_{\mathcal C}.
\]
\end{La}

The next theorem characterizes the reproducing kernel Hilbert
space associated with $K_r(t,s)$. It is a direct consequence of
formula \eqref{marseille} and Lemma \ref{Vaugirard}.

\begin{Tm}
\label{totoche}
The reproducing kernel Hilbert space ${\mathcal H}(K_r)$
associated with $K_r(t,s)$ consists of functions of the form
\begin{equation}
{\mathcal H}(K_r)=\left\{\, F\,:\, F(t)=\int_{\mathbb
R}\dfrac{e^{itu}-1}{u}f(u)dm(u),\right\}
\label{St-Michel}
\end{equation}
where $f$ is in the closed linear span in ${\mathbf L}_2(dm)$ of
the functions
\begin{equation}
\chi_s(u)=\frac{e^{isu}-1}{u},\quad s\in{\mathbb R},
\label{haifa}
\end{equation}
with the norm
\begin{equation}
\label{napoleon} \|F\|_{{\mathcal H}(K_r)}=\|f\|_{{\mathbf
L}_2(dm)}.
\end{equation}
\label{attention-les-yeux}
\end{Tm}

Theorem \ref{totoche} is proved through different methods
in \cite[Theorem 4.1 p.
948]{fbm88} for $dm(u)=\frac{1}{\pi}|u|^{1-2H}du$ with
$H\in(1/2,1)$. See also \cite{MR1280932},
\cite{MR2145497}. We note that in the case $H\in(1/2,1)$ the
covariance kernel can also be expressed in terms of double
integrals as follows: one has
\[
|t|^{2H}+|s|^{2H}-|t-s|^{2H}=2H(2H-1)\int_0^t\int_0^s
|u-v|^{2H-2}dudv
\]
and the associated reproducing kernel Hilbert space can be
expressed in terms of the preHilbert space of functions $f$ such
that
\[
\iint_{{\mathbb R}_+^2}f(u)f(v)|u-v|^{2H-2}dudv<\infty,
\]
See \cite{MR1801485}, \cite{duncan1}.
This last space is not complete; see \cite[p. 270]{MR1790083}.\\

{\bf Remarks:}\\
$(a)$ In \eqref{St-Michel}, one can take $f\in{\mathbf L}_2(m)$,
rather than in the stated closed linear span; one then looses in
general the uniqueness of $f$ in the representation of $F$ and one
has to change \eqref{napoleon} to $\inf \|f\|_{{\mathbf
L}_2(m)}$, where the $\inf$ is over all $f$ corresponding to the
given function
$F$.\\
$(b)$ An important case is when $dm$ is absolutely continuous with
respect to Lebesgue measure and when moreover its derivative
satisfies
\[
\int_{\mathbb R}\dfrac{\ln m^\prime(u)}{1+u^2}du>-\infty
\]
Then $m^\prime$ admits an outer factorization $m^\prime=|h|^2$;
see \cite{Dmk}.\\
$(c)$ For any $m$ such that the space
${\mathbf L}_2(dm)$ is infinite dimensional and any $T>0$, the
closed linear span of the functions $\chi_t$ (defined by
\eqref{haifa}) for $|t|\le T$ has a special structure. It is a
reproducing kernel space of {\it entire functions} of the kind
introduced by L. de Branges; see \cite{MR38:6386}, \cite{dbbook},
\cite{dym-70}, \cite{Dmk}. Its reproducing kernel is of the form
\[
\frac{a(z)a(w)^*-b(z)b(w)^*}{-i(z-w^*)}
\]
where the functions $a$ and $b$ are entire functions of finite
exponential type. In the case of the fractional Brownian motion,
these are the homogeneous de Branges spaces of entire functions
(see \cite[pp. 184-189]{dbbook}). The functions $a$ and $b$ are
then expressed in terms of Bessel functions, a fact already
observed (in a slightly
different language) by Krein in \cite{krein_1954}. See also
\cite{MR2145497}. We refer to \cite{ad1} for more information
on the reproducing kernel spaces of the kind introduced by
de Branges and Rovnyak.\\
$(d)$ To prove that the kernels $K_{H,1}$ are positive, one can
also use the theory of negative kernels (see \cite{MR86b:43001}).
Another proof can be found in \cite[p. 209]{Mallat-francais}. The
associated reproducing kernel Hilbert space is given in for
instance in \cite{MR0242247} (with different representations of the kernels
corresponding to  $H<1/2$ and $H>1/2$; see
\cite[p. 136]{MR0242247}), in \cite[p. 319]{berlinet}, and can
also be deduced from the model for the fractional Brownian motion
given in e.g. \cite[p. 50]{Lifshits95}, \cite[\S 7.2 p.
318]{MR1280932}, \cite{benassi_cohen}. A characterization as
operator range can be deduced from \cite[Lemma 3.1 p. 182]{DS}).
When $H=1/2$, it is well known that the associated reproducing
kernel Hilbert space is the Sobolev space of functions absolutely
continuous on the real line and with square summable derivatives;
see \cite[p.
25]{Malliavin}, \cite{MR99f:60082}.\\

We now study a relation between the functions
$f$ and $F$ of \eqref{St-Michel}. We restrict ourselves to the case
where $m$ is
absolutely continuous with respect to Lebesgue measure, and
resort to the theory of distributions, as we now explain. Denote
by ${\mathcal S}$ the space of rapidly decreasing functions. These
are the infinitely differentiable functions $\sigma(u)$ such that
for every choice of $(k,r)\in{\mathbb N}^2$,
\[
\sup_{u\in{\mathbb R}}|(1+u^2)^r\sigma^{(k)}(u)|<\infty.
\]
We refer the reader to \cite{Schwartz66} and \cite{Treves67} for
more information on this space and in particular on its topology.
We denote by ${\mathcal S}^\prime$ the set of continuous linear
functionals on ${\mathcal S}$. The elements of ${\mathcal
S}^\prime$ are called {\it tempered distributions}.

\begin{Tm}
A function $F$ of the form \eqref{St-Michel}, defines an element
in ${\mathcal S}^\prime$ whose derivative is the Fourier
transform of the distribution
\begin{equation}
\label{Chevaleret} T_{m^\prime f}\sigma=\int_{\mathbb
R}\sigma(u)m^\prime(u)f(u)du.
\end{equation}
\end{Tm}

{\bf Proof:} We first note that the function $m^\prime f$ indeed
defines an element in ${\mathcal S}^\prime$ via the formula
\eqref{Chevaleret} since, for every $\sigma\in{\mathcal S}$, the integral
\[
\begin{split}
\int_{\mathbb R}m^\prime(u)f(u){\sigma}(u)du&=\int_{\mathbb
R}\dfrac{m^\prime(u)}{(1+u^2)^\ell}f(u)(1+u^2)^\ell{\sigma}(u)du,
\end{split}
\]
exists, 
where $\ell\in{\mathbb N}$ is such that $(1+u^2)^\ell \sigma(u)$ is
bounded on the real line.
Using Fubini's theorem and integration by parts, we have for
${\sigma}\in{\mathcal S}$
\[
\begin{split}
D_F{\sigma}&=\int_{\mathbb R}{\sigma}^{\prime}(t)F(t)dt\\
&=\int_{\mathbb R}{\sigma}^\prime(t)\left(\int_{\mathbb
R}\dfrac{e^{itu}-1}{u}m^\prime(u)f(u)du\right)dt\\
&=\int_{\mathbb R}\left(\int_{\mathbb
R}{\sigma}^\prime(t)\dfrac{e^{itu}-1}{u}dt\right)m^\prime(u)f(u)du\\
&=-\int_{\mathbb R}\widetilde{{\sigma}}(u)m^\prime(u)f(u)du,
\end{split}
\]
where we denote by $\widetilde{\sigma}$ the Fourier transform
\[
\widetilde{\sigma}(w)=\int_{\mathbb R}\sigma(t)e^{itw}dt.
\]

This last equality expresses the fact that
$\widetilde{D_F}=T_{m^\prime f}$, and that the distributional
derivative of $F$ is the distributional Fourier transform of
$m^\prime f$.
\mbox{}\qed\mbox{}\\

Let us now restrict the above results to $dm(u)=du$. Then, we note
that
\[
\int_0^\infty\dfrac{1-\cos(tu)}{u^2}\,du=\dfrac{\pi|t|}{2}
\]
(as is seen by integrating the function
$f(z)=\frac{1-e^{itz}}{z^2}$ around an appropriate contour), and
so
\[
\begin{split}
\dfrac{1}{2}\left\{|t|+|s|-|t-s|\right\}&=\dfrac{1}{\pi}\int_0^\infty
\dfrac{\cos((t-s)u)-\cos(tu)-\cos(su)+1}{u^2}\,du\\
&={\rm Re}~\dfrac{1}{2\pi}\int_{\mathbb
R}\dfrac{e^{i(t-s)u}-e^{itu}-e^{-isu}+1}{u^2}du\\
&=\dfrac{1}{2\pi}\int_{\mathbb
R}\dfrac{e^{i(t-s)u}-e^{itu}-e^{-isu}+1}{u^2}du\\
&=\dfrac{1}{2\pi}\int_{\mathbb R}
\dfrac{e^{itu}-1}{u}\dfrac{e^{-isu}-1}{u}\,du\\
&=\frac{1}{2\pi}\langle \chi_t,\chi_s\rangle_{{\mathbf L}_2(du)}.
\end{split}
\]
The closed linear span in ${\mathbf L}_2(du)$ of the functions
$\chi_t$ ($t\in{\mathbb R}$) is equal to ${\mathbf L}_2(du)$, and
the reproducing kernel Hilbert space associated with the kernel
$\frac{1}{2}\left\{|t|+|s|-|t-s|\right\}$ is equal to the set of
functions $F$ of the form
\[
F(t)=\int_{\mathbb R}\dfrac{e^{itu}-1}{u}f(u)du,
\]
where $f\in{\mathbf L}_2(du)$, with the norm
$\|F\|=\frac{1}{\sqrt{2\pi}}
\|f\|_{{\mathbf L}_2(du)}$. This space is known to be the
Sobolev space. We regenerate this result as follows. For
$f\in{\mathbf L}_1(du)\cap{\mathbf L}_2(du)$, Fubini's theorem
shows that
\begin{equation}
\label{2.5}
F(t)=\int_{\mathbb R}\dfrac{e^{itu}-1}{u}f(u)du
=\int_0^t \left(\int_{\mathbb R}e^{ius}f(u)du\right)ds
=\int_0^t\widetilde{f}(s)ds,
\end{equation}
and
\[
\|F\|=\frac{1}{\sqrt{2\pi}}\|f\|_{{\mathbf L}_2({\mathbb R})}=
\|\widetilde{f}\|_{{\mathbf L}_2({\mathbb R})}.
\]
The result for general $f\in{\mathbf L}_2(du)$
follows by approximation. If we restrict $t,s\in[0,T]$ for some
$T>0$, then $f$ belongs to the closed linear span in ${\mathbf
  L}_2({\mathbb R})$ of the functions $\chi_s$, $s\in[-t,0]$. Thus
$f(t)=\int_0^T e^{-iut}g(u)du$ for some
$g\in{\mathbf L}_2([0,T])$,
and we have
\[F(t)=\int_0^t\left(\int_{\mathbb
  R}e^{isx}\left(\int_0^Te^{-ixu}g(u)du\right)
dx\right)ds=\int_0^tg(u)du.\]

\section{Integrals and products of positive functions}
\label{clementine}
\setcounter{equation}{0}
In this section we briefly review a number of results related to
reproducing kernel Hilbert spaces. The results are stated for
functions positive on the real line because this is the setting
we consider, but are valid for functions positive on any set. The
first result is in Schwartz's paper \cite[p. 170]{schwartz}. We
outline the proof in our special case for completeness.

\begin{Tm}
Let $k_u(t,s)=k(t,s,u)$ be a continuous function on ${\mathbb
R}^2\times{\mathbb R}_+$ such that for every $u\in{\mathbb R}_+$
the kernel $k(t,s,u)$ is positive on ${\mathbb R}$, and let $q$ be a
function positive and increasing on ${\mathbb R}_+$. Assume that
$\int_0^\infty k(t,t,u)dq(u)<\infty$.
Then the function 
\begin{equation}
\mathbf k(t,s)=\int_0^\infty k(t,s,u)dq(u)
\label{Temple, ligne 3}
\end{equation}
is positive on the real line and
its reproducing kernel Hilbert space
is the Hilbert space $\int_0^\infty{\mathcal
H}(k_u)dq(u)$. 
\end{Tm}
{\bf Proof:} First we note that the function
\eqref{Temple, ligne 3} is well defined. Indeed, from the
positivity of $k_u$ we get that
\[
|k(t,s,u)|^2\le k(t,t,u)k(s,s,u)
\]
and hence by Cauchy-Schwarz inequality we have
\[
\begin{split}
|\int_0^\infty k(t,s,u)dq(u)|&\le\int_0^\infty
k(t,t,u)^{1/2}k(s,s,u)^{1/2}dq(u)\\
&\le\left(\int_0^\infty
k(t,t,u)dq(u)\right)^{1/2}\left(\int_0^\infty
k(s,s,u)dq(u)\right)^{1/2}<\infty.
\end{split}
\]
Consider now the vector space ${\mathcal M}(\mathbf k)$ spanned by the
finite linear combinations of the functions ${\mathbf k}(t,s)$,
with inner product
\[
\langle {\mathbf k}(\cdot, w),{\mathbf k}(\cdot, v)
\rangle_{{\mathcal M}({\mathbf k})}={\mathbf k}(v,w).
\]
An element $F\in{\mathcal M}({\mathbf k})$ can thus be written as
a finite sum
\[
F(t)=\sum_{\ell=1}^m \mathbf k(t,s_\ell)c_\ell=\int_0^\infty
\left(\sum_{\ell=1}^m k(t,s_\ell,u)c_\ell\right)dq(u),\quad
c_\ell\in{\mathbb C},\]
with norm
\[
\|F\|^2_{{\mathcal M}({\mathbf k})}=\int_0^\infty\|
\sum_{\ell=1}^m k(\cdot,s_\ell,u)c_\ell\|^2_{{\mathcal
H}(k_u)}dq(u).
\]
The space ${\mathcal H}({\mathbf k})$ is on the one hand the
completion of ${\mathcal M}({\mathbf k})$ with respect to this
inner product. By \cite[\S5, p. 146]{MR0352996}, this is also the
description of the space $\int_0^\infty {\mathcal H}(k_u)dq(u)$.
By uniqueness of the completion we have
the result. \qed\mbox{}\\

Elements of the space $\int_0^\infty {\mathcal H}(k_u)dq(u)$ can be
described as follows (see \cite[p. 170]{schwartz}, \cite[pp. 77 and
96]{saitoh}): A function $F$ belongs to ${\mathcal
H}({\mathbf k})$ if and only if it can be written as
\begin{equation}
\label{La_Muette} F(t)=\int_0^\infty x(t,u)dq(u)
\end{equation}
where for every $u$, the function $x_u\,:\,t\mapsto
x(t,u)\in{\mathcal H}(k_u)$ and the map $u\mapsto x_u$ is weakly
measurable. Moreover,
\begin{equation}
\label{perlette} \|F\|_{{\mathcal H}({\mathbf
k})}^2=\inf\int_0^\infty \|x_u\|^2_{{\mathcal H}(k_u)}dq(u)
\end{equation}
where the infimum is over all the representations of a given $F$
in \eqref{La_Muette}.\\

We refer the reader to \cite{MR0044536}, \cite{MR0056611},
\cite{MR0029101} and \cite[pp. 139--161]{MR0352996} for more
information on integrals of Hilbert spaces.\\

As a special case we have:

\begin{Tm}
\label{Chevaleret, ligne 6} Let $k_n(t,s)$, $n=1,2,\ldots$ be
positive functions on ${\mathbb R}$ and let
\[
\mathbf k(t,s)=\sum_{n=1}^\infty k_n(t,s).
\]
Assume that $\mathbf k(t,t)<\infty$ for all $t\in{\mathbb R}$.
 Then
\begin{equation}
\label{helene}
{\mathcal H}(\mathbf k)=\sum_{n=1}^\infty{\mathcal H}(k_n)
\end{equation}
in the following sense: ${\mathcal H}(\mathbf k)$ consists of
the functions which can be written as
\begin{equation}
\label{Rue du Bac} f(t)=\sum_{n=1}^\infty f_n(t),
\end{equation}
where $f_n\in{\mathcal H}(k_n)$ and
\begin{equation}
\sum_{n=1}^\infty\|f_n\|^2_{{\mathcal H}(k_n)}<\infty.
\label{Rue_du_Bac2}
\end{equation}
Moreover, we have:
\begin{equation}
\label{quartier_latin} 
\|f\|_{{\mathcal H}(\mathbf
k)}^2=\inf\sum_{n=1}^\infty\|f_n\|^2_{{\mathcal H}(k_n)},
\end{equation}
where the $\inf$ is over all the possible decompositions
\eqref{Rue du Bac} satisfying \eqref{Rue_du_Bac2}.
\end{Tm}

We note that the infimum is in fact a minimum and is achieved for
the choice
\[
f_n=i_ni_n^*f
\]
where $i_n$ is the inclusion map $f\mapsto f$ from ${\mathcal
H}(k_n)$ into ${\mathcal H}(\mathbf k)$. See \cite{aron}, \cite{saitoh}.\\

The sum \eqref{helene} is in general not direct, and the spaces
are called complementary; see \cite{dbjfa} for more information on
complementation theory.\\

The following auxiliary result will be instrumental in the
following sections. Note that the result is stated in a form
convenient for our present purpose.

\begin{Tm}
\label{rosemary}
 Let $k_n(t,s,u), n=1,2,\ldots$ be a sequence of continuous functions
defined on ${\mathbb R}^2\times {\mathbb R}_+$ and assume that
for every $u\in{\mathbb R}_+$ and every positive integer $n$ the
function $k_n(t,s,u)$ is positive on ${\mathbb R}$. Assume that
\[
\sum_{n=1}^\infty k_n(t,s,u)
\]
is continuous on ${\mathbb R}^2\times {\mathbb R}_+$ and that for
every $t\in{\mathbb R}$
\[
\int_0^\infty\left(\sum_{n=1}^\infty k_n(t,t,u)\right)du<\infty
\]
Then,
\[
\int_0^\infty \left(\sum_{n=1}^\infty k_n(t,s,u)\right)du=
\sum_{n=1}^\infty\int_0^\infty k_n(t,s,u)du.
\]
\end{Tm}
{\bf Proof:} For every $u\in{\mathbb R}_+$ and $n\in{\mathbb
N}^*$ we have
\[
|k_n(t,s,u)|^2\le k_n(t,t,u)k_n(s,s,u),\]
and so,
\[
\begin{split}
|\sum_{n=1}^\infty k_n(t,s,u)|&\le\sum_{n=1}^\infty
 \sqrt{k_n(t,t,u)}\sqrt{k_n(s,s,u)}\\
 &\le\left(\sum_{n=1}^\infty k_n(t,t,u)\right)^{1/2}\left(\sum_{n=1}^\infty
 k_n(s,s,u)\right)^{1/2}
\end{split}
\]
and the function $u\mapsto \sum_{n=1}^\infty k_n(t,s,u)$ is absolutely
summable for every choice of $t,s\in{\mathbb R}$.
The result is then a direct consequence of the dominated
convergence theorem.\mbox{}\qed\mbox{}\\

The next theorem is classical; see for instance
\cite{aron},\cite{saitoh}.
\begin{Tm}
Let $k_1$ and $k_2$ be two functions positive on the real line.
Then the function $k_1k_2$ is also positive on the real line.
Assume moreover that the functions $k_1$ and $k_2$ are continuous.
Then the Hilbert spaces ${\mathcal H}(k_1)$ and ${\mathcal
H}(k_2)$ are separable. The space ${\mathcal H}(k_1k_2)$ consists
of the restrictions on the diagonal of the elements of tensor
product ${\mathcal H}(k_1)\otimes{\mathcal H}(k_2)$. Let
$(f_n)_{n=0,\ldots}$ and $(g_n)_{n=0,\ldots}$ be Hilbert space
basis of the spaces ${\mathcal H}(k_1)$ and ${\mathcal H}(k_2)$
respectively. Then ${\mathcal H}(k_1k_2)$ consists of the
functions of the form
\[
F(t)=\sum_{n,m=0}^\infty c_{n,m}f_n(t)g_m(t)\]
with
\[
\|F\|^2_{{\mathcal H}(k_1k_2)}=\sum_{n,m=0}^\infty |c_{n,m}|^2.
\]
\end{Tm}
From this result we obtain that $f_1f_2\in{\mathcal H}(k_1k_2)$ where
$f_\ell\in{\mathcal H}(k_\ell)$, $\ell=1,2$, and
\begin{equation}
\label{saitoh_1}
\|f_1f_2\|_{{\mathcal H}(k_1k_2)}\le \|f_1\|_{{\mathcal
H}(k_1)}\|f_2\|_{{\mathcal H}(k_2)}\quad\mbox{{\rm and that}}\quad
\|f^n\|^2_{{\mathcal H}(k^n)}\le \|f\|^{2n}_{{\mathcal H}(k)}
\end{equation}
where $k$ is a positive kernel. See \cite[p 243]{saitoh2} for the
second inequality.
\section{Positivity of the kernels $K_{H,1}(t,s)$}
\label{Avenue Foch}
\setcounter{equation}{0}
In  this section we present two decompositions 
of the form \eqref{La_Fourche} of the kernels
$K_{H,1}(t,s)$, defined in \eqref{ioanna}. The first is valid
for any $H\in(0,1)$ while the second is valid only for
$H\in(0,1/2)$.

\begin{Tm} Let $0<H<1$. \\
$(a)$ It holds that
\begin{equation}
\label{total-recall-nice-1991}
\begin{split}
\frac{\Gamma(1-H)}{2H}
(|t|^{2H}+|s|^{2H}-|t-s|^{2H})&=\\
&\hspace{-2.5cm} =\int_0^\infty
\frac{(1-e^{-u^2t^2})(1-e^{-u^2s^2})}{u^{1+{2H}}}du+\\
&\hspace{-2cm}+
 \sum_{n=1}^\infty
\frac{2^{n-1}\Gamma(n-H)}{n!}\frac{t^ns^n}{(t^2+s^2)^{n-H}}.
\end{split}
\end{equation}
$(b)$A function $F$ belongs to ${\mathcal
H}(\frac{2H}{\Gamma(1-H)}K_{H,1})$ if and only if it can be
written as
\[
F(t)=\sum_{n=0}^\infty F_n(t)
\]
with
\[
F_0(t)=\int_0^\infty \frac{(1-e^{-u^2t^2})x_0(u)}{u^{1+2H}}du,
\]
where $x_0$ belongs to the closed linear span of the functions
$1-e^{-u^2s^2}$ ($s\in{\mathbb R}$) in ${\mathbf
L}_2(1/u^{1+{2H}})$, and where for every $n\ge 0$,
\[
F_n(t)=t^n\int_0^\infty e^{-u^2t^2}x_n(u)u^{2n-1-2H}du,\] where
$x_n$ belongs to the closed linear span of the functions
$e^{-u^2s^2}$ ($s\in{\mathbb R}$) in ${\mathbf
L}_2(u^{2n-1-{2H}})$. For any such decomposition,
\begin{equation}
\label{Porte d'Orleans, ligne 4}
\|F\|^2\le\sum_{n=0}^\infty\frac{n!}{2^n}\|x_n\|^2_{{\mathbf
L}^2(u^{2n-1-2H})},
\end{equation}
and there is a unique decomposition for which equality hold in
\eqref{Porte d'Orleans, ligne 4}.
\end{Tm}
{\bf Proof:} The proof is based on an idea of I. Schoenberg. We
use the equality (see \cite[equation (8) p. 526]{MR1501980})
\begin{equation}
 |t|^{2H}=\dfrac{\int_0^\infty (1-e^{-u^2t^2})u^{-1-{2H}}du}
{\int_0^\infty (1-e^{-u^2})u^{-1-{2H}}du}. \label{schoenberg}
\end{equation}
The equality itself is proved by using the change of variable
$u\mapsto u|t|$ in the integral
\[\int_0^\infty
(1-e^{-u^2t^2})u^{-1-{2H}}du,\quad 0<H<1.\]
See \cite[p. 526]{MR1501980}. Moreover, integration by part leads
to
\[
\int_0^\infty (1-e^{-u^2})u^{-1-{2H}}du=\dfrac{\Gamma(1-H)}{2H}.
\]
From \eqref{schoenberg} it follows that
\begin{equation}
\label{total-recall}
\begin{split}
\left(\frac{\Gamma(1-H)}{2H}\right)
(|t|^{2H}+|s|^{2H}-|t-s|^{2H})&=\\
&\hspace{-2.5cm} =\int_0^\infty
\frac{(1-e^{-u^2t^2})(1-e^{-u^2s^2})}{u^{1+{2H}}}du+\\
&\hspace{-2cm}+
\int_0^\infty\frac{e^{-u^2t^2}(e^{2u^2ts}-1)e^{-u^2s^2}}{u^{1+{2H}}}du.
\end{split}
\end{equation}
The second term is a decomposition of the form \eqref{Temple,
ligne 3} with
\[
k(t,s,u)=e^{-u^2t^2}(e^{2u^2ts}-1)e^{-u^2s^2}\quad{\rm and}\quad
q^\prime(u)=\frac{1}{u^{1+{2H}}}.
\]
The positivity of $K_{H,1}(t,s)$ follows by using the fact that
sums and products of positive functions are still positive. We
note that equation \eqref{total-recall} proves that
$K_{H,1}(t,s)$ is a covariance function, but it does not provide
a {\it model} for it (in the sense given in \cite[p.
41]{Lifshits95}). Still, one can be somewhat more explicit by
using formulas \eqref{La_Muette} and \eqref{perlette}: A function
$F$ belongs to the reproducing kernel Hilbert space with
reproducing kernel
\[
M(t,s)=
\int_0^\infty\frac{e^{-u^2t^2}(e^{2u^2ts}-1)e^{-u^2s^2}}{u^{1+{2H}}}du
\]
if and only if it can be written as
\[
F(t)=\int_0^\infty \frac{f(t,u)}{u^{1+2H}}du
\]
where for every $u$ the function $f(t,u)$ belongs to the
reproducing kernel Hilbert space ${\mathcal H}(N_u)$ with reproducing
kernel
\[
N_u(t,s)= e^{-u^2t^2}(e^{2u^2ts}-1)e^{-u^2s^2}
\]
and
\[
\int_0^\infty\frac{\|f(t,\cdot)\|^2_{{{\mathcal H}(N_u)}}}{{u^{1+{2H}}}}
du<\infty.
\]
Next, and using Theorem \ref{rosemary}, write:
\[
\begin{split}
\int_0^\infty\frac{e^{-u^2t^2}(e^{2u^2ts}-1)e^{-u^2s^2}}{u^{1+{2H}}}du&=
\sum_{n=1}^\infty \frac{2^nt^ns^n}{n!}\int_0^\infty
e^{-u^2(t^2+s^2)}u^{2n-1-2H}du\\
&= \sum_{n=1}^\infty
\frac{2^{n-1}\Gamma(n-H)}{n!}\frac{t^ns^n}{(t^2+s^2)^{n-H}},
\end{split}
\]
which recovers \eqref{total-recall-nice-1991}.
\mbox{}\qed\mbox{}\\

For instance, for a given real number $s$ and for
\[
F(t)=\frac{2H}{\Gamma(1-H)}K_{H,1}(t,s)
\]
we have \( x_0(u)=(1-e^{-u^2s^2})\) and for $n\ge 1$,
\[
x_n(u)=\frac{2^n}{n!}s^ne^{-u^2s^2}.\]
For this choice of $x_0,x_1,\ldots$ equality holds in
\eqref{Porte d'Orleans, ligne 4}.\\

The second approach is also related to Schoenberg's work. Let
$dm(u)$ be a positive measure on $[0,\infty)$ subject to the
condition
\[
\int_1^\infty dm(u)<\infty.
\]
(Note that the Lebesgue measure does not satisfy this last
requirement.)\\

The function
\begin{equation}
\label{new_r}
r(t)=\int_0^\infty (1-e^{-u|t|})dm(u)
\end{equation}
is then defined for all real $t$'s. Such functions $r(t)$ have
been introduced in \cite[Definition 2, p. 825]{MR1503439}.
Functions of the form $\sqrt{r(t^2)}$ are exactly the functions
such that a separable real Hilbert space ${\mathcal H}$ (with
norm denoted by $\|\cdot\|$), endowed with the metric $F(\|x\|)$,
can be isometrically imbedded in ${\mathcal H}$. See
\cite[Theorem 6, p. 828]{MR1503439}. In particular, $r$ is of the
form \eqref{eureka}.

\begin{Tm}
Let $r(t)$ be of the form \eqref{new_r}. Then the kernel
$K_r(t,s)$:
\[
K_r(t,s)=r(t)+r(s)-r(t-s).
\]
is positive on the real line. Moreover, with
$K_{1,1}(t,s)=|t|+|s|-|t-s| \stackrel{\rm def.}{=}K(t,s)$,
\begin{equation}
\label{St Philippe du Roule, ligne 9}
\begin{split}
K_r(t,s) &=\int_0^\infty(1-e^{-u|t|})(1-e^{-u|s|})dm(u)+\\
&\hspace{5mm}+\sum_{n=1}^\infty\dfrac{K^n(t,s)}{n!}
\left(\int_0^\infty u^ne^{-u(|t|+|s|)}dm(u)\right).
\end{split}
\end{equation}
\end{Tm}
{\bf Proof:} Using Theorem \ref{rosemary} we have:
\[
\begin{split}
K_r(t,s)&=\int_0^\infty \left\{
1-e^{-u|t|}+1-e^{-u|s|}+e^{-u|t-s|}-1\right\}
m(u) du\\
&=\int_0^\infty(1-e^{-u|t|})(1-e^{-u|s|})dm(u)+ \int_0^\infty
\left\{e^{-u|t-s|}-e^{-u(|t|+|s|)}\right\}dm(u)\\
&=\int_0^\infty(1-e^{-u|t|})(1-e^{-u|s|})dm(u)+\\
&\hspace{5mm}+ \int_0^\infty e^{-u|t|}
(e^{u(|t|+|s|-|t-s|)}-1)e^{-u|s|}dm(u).
\end{split}
\]
which proves the positivity of the kernel since for every $u\ge
0$, one has
\[
e^{u(|t|+|s|-|t-s|)}-1=\sum_{1}^\infty
\dfrac{u^n}{n!}(|t|+|s|-|t-s|)^n
\]

\mbox{}\qed\\

{\bf Remarks:}\\
$(a)$ The function
\begin{equation}
\label{segolene} \Gamma(z,\mu)=\int_0^\infty u^{z-1}e^{-\mu
u}dm(u)
\end{equation}
is the generalization for any $m(u)$ of the classical Gamma
function. Furthermore, the decomposition \eqref{St Philippe du
Roule, ligne 9}
is conducive in applications with regard to non--linear
transforms. See Section \ref{NL}.\\
$(b)$ The formula \eqref{St Philippe du Roule, ligne 9} allows to
characterize the associated reproducing kernel Hilbert space.
\\
$c)$ In view of the result of von Neumann and Schoenberg
mentioned in the introduction, and as we already noticed above,
any function of the form
\eqref{new_r} is of the form \eqref{eureka}.\\

In the case where
\[
dm(u)=\dfrac{e^{-u}}{u}du
\]
we have $r(t)=\ln(1+|t|)$. This follows from the formula in
\cite[Corollary 2.10]{MR86b:43001}:
\[
\ln (1+z)=\int_0^\infty(1-e^{-uz})\dfrac{e^{-u}du}{u},\quad {\rm
Re}~z\ge 0.
\]
It follows that the function
\[
q(t,s)=\ln(1+|t|)+\ln(1+|s|)-\ln(1+|t-s|)
\]
is positive on ${\mathbb R}$.\\

The choice
\[
dm(u)=\dfrac{2H}{\Gamma(1-2H)} \dfrac{1}{u^{1+2H}}du
\]
corresponds to the covariance function of the
fractional Brownian motion for $H\in(0,1/2)$
since (see \cite[Corollary 2.10 p. 78]{MR86b:43001}),
\[
z^{2H}=\dfrac{2H}{\Gamma(1-2H)}\int_0^\infty
(1-e^{-uz})\dfrac{du}{u^{1+2H}},\]
where $H\in(0,1/2)$ and ${\rm Re}~z\ge 0$. We check this formula
for positive $t$ as follows. The change of variable $v=tu$ leads
to
\[
\int_0^\infty (1-e^{-uz})\dfrac{du}{u^{1+2H}}=t^{2H}\int_0^\infty
(1-e^{-u})\dfrac{du}{u^{1+2H}}=\frac{\Gamma(1-2H)}{2H}
\]
where the last equality is obtained by integration by parts.\\

Specializing \eqref{St Philippe du Roule, ligne 9} to this case
leads to:

\begin{Tm} Let $H\in(0,1/2)$. Then:
\begin{equation}
\label{ohohoh}
\begin{split}
|t|^{2H}+|s|^{2H}-|t-s|^{2H}
&=2H(|t|+|s|-|t-s|)+\\
&\hspace{5mm}+\frac{2H}{\Gamma(1-2H)}\left\{
\int_0^\infty\frac{(1-e^{-u|t|})(1-e^{-u|s|})}{u^{1+2H}}du+\right.\\
&\hspace{3cm}+\left.
\sum_{n=2}^\infty\dfrac{K^n(t,s)\Gamma(n-2H)}{n!(|t|+|s|)^{n-2H}}\right\}.
\end{split}
\end{equation}
In particular, the space ${\mathcal H}(K_{H,1})$ contains
contractively a copy of ${\mathcal H}(K_{1,1})$.
\label{surprise}
\end{Tm}
{\bf Proof:} Indeed, using Theorem \ref{rosemary} we have:
\begin{equation*}
\begin{split}
|t|^{2H}+|s|^{2H}-|t-s|^{2H}&=\frac{2H}{\Gamma(1-2H)}\left\{
\int_0^\infty\frac{(1-e^{-u|t|})(1-e^{-u|s|})}{u^{1+2H}}du+\right.\\
&\hspace{3cm}+\left.
\sum_{n=1}^\infty\dfrac{K^n(t,s)\Gamma(n-2H)}{n!(|t|+|s|)^{n-2H}}
\right\}\\
&=2H(|t|+s|-|t-s|)+\\
&\hspace{5mm}+\frac{2H}{\Gamma(1-2H)}\left\{
\int_0^\infty\frac{(1-e^{-u|t|})(1-e^{-u|s|})}{u^{1+2H}}du+\right.\\
&\hspace{3cm}+\left.
\sum_{n=2}^\infty\dfrac{K^n(t,s)\Gamma(n-2H)}{n!(|t|+|s|)^{n-2H}}\right\}.
\end{split}
\end{equation*}
\mbox{}\qed\mbox{}\\

{\bf Remark:} Formula \eqref{total-recall-nice-1991} is valid for
any $0<H<1$ and we see that terms of the form $\Gamma(n-H)$
appear in it ($n\ge 1$). On the other hand, formula
\eqref{ohohoh} is valid only for $0<H<1/2$, with terms of the form
$\Gamma(n-2H)$.

\section{A general family of positive functions}
\setcounter{equation}{0}
\begin{Tm}
Let $r(t)$ be a real valued function such that the kernel
$K_r(t,s)$, defined by \eqref{mathilde_est_revenue}, is positive
on the real line, and let $\varphi$ is of the form \eqref{new_r}.
Then:
\begin{equation}
\label{jezabel} \varphi(r(t)+r(s))-\varphi(r(t-s)) =
\sum_{n=1}^\infty\dfrac{K_r(t,s)^n}{n!}\Gamma(n, r(t)+r(s))
\end{equation}
where $\Gamma(z,\mu)$ has been defined in \eqref{segolene}.
In particular the function
\begin{equation}
\label{Bercy} 
V(t,s)=\varphi(r(t)+r(s))-\varphi(r(t-s))
\end{equation}
is positive on the real line.
\end{Tm}

{\bf Proof:} With $\varphi$ of the form \eqref{new_r}, one has

\begin{equation*}
\begin{split}
\varphi(r(t)+r(s))-\varphi(r(t-s))&=\int_0^\infty \left\{
e^{-ur(t-s)}-e^{-u(r(t)+r(s))}\right\}dm(u)\\
&=\int_0^\infty e^{-ur(t)})\left\{
e^{u(r(t)+r(s)-r(t-s))}-1\right\}e^{-ur(s)}dm(u)\\
&= \sum_{n=1}^\infty\dfrac{K_r(t,s)^n}{n!}\Gamma(n, r(t)+r(s)).
\end{split}
\end{equation*}
This allows to conclude the proof since the functions
$K_r(t,s)^n$ and $\Gamma(n, r(t)+r(s))$ are positive on the real
line, and so is
their product. \mbox{}\qed\mbox{}\\

We now consider the case $\varphi(z)=z^\alpha$ and obtain the
following analogue of Theorem \ref{surprise}:

\begin{Tm} Let $\alpha\in(0,1)$. Then:
\begin{equation}
\label{Monceau, ligne 2}
\begin{split}
(r(t)+r(s))^\alpha-r(t-s)^\alpha&= \alpha
K_r(t,s)+\frac{\alpha}{\Gamma(1-\alpha)}
\sum_{n=2}^\infty\frac{K_r(t,s)^n\Gamma(n-\alpha)}{n!(r(t)+r(s))^{n-\alpha}}.
\end{split}
\end{equation}
In particular, the space with reproducing kernel
$(r(t)+r(s))^\alpha-r(t-s)^\alpha$ contains contractively a copy
of ${\mathcal H}(K_r)$.
\end{Tm}
{\bf Proof:} We have already recalled that
\[
z^\alpha=\frac{\alpha}{\Gamma(1-\alpha)}\int_0^\infty(1-e^{-uz})
\frac{du}{u^{1+\alpha}}.
\]
Using the preceding arguments we therefore obtain:

\begin{equation*}
\begin{split}
(r(t)+r(s))^\alpha-r(t-s)^\alpha&=
\frac{\alpha}{\Gamma(1-\alpha)}\sum_{n=1}^\infty\frac{K_r(t,s)^n}{n!}
\int_0^\infty u^ne^{-u(r(t)+r(s))}\frac{du}{u^{1+\alpha}}\\
&=
\frac{\alpha}{\Gamma(1-\alpha)}\sum_{n=1}^\infty\frac{K_r(t,s)^n}{n!}
\dfrac{1}{(r(t)+r(s))^{n-\alpha}}\Gamma(n-\alpha)\\
&=\alpha K_r(t,s)+\frac{\alpha}{\Gamma(1-\alpha)}
\sum_{n=2}^\infty\frac{K_r(t,s)^n\Gamma(n-\alpha)}{n!(r(t)+r(s))^{n-\alpha}}.
\end{split}
\end{equation*}
The last claim of the theorem follows from
\eqref{Monceau, ligne 2} and Theorem \ref{Chevaleret, ligne 6}
\mbox{}\qed\mbox{}\\

{\bf Remark:} Formula \eqref{Monceau, ligne 2} is still valid for
$\alpha=1$ since the function $\Gamma(z)$ has a pole at the
origin.

\section{Associated non-linear transforms}
\setcounter{equation}{0}
\label{NL}

We now associate a non-linear transform (in the sense of Saitoh;
see \cite[Appendix 2, p. 243] {saitoh2}) to the decomposition
\eqref{jezabel}. We begin with a general result.

\begin{Tm}
Let $K$ and $K_n,\, n=1,2,\ldots$ be positive
kernels on the real line and define
\[
{\mathbf K}(t,s)=\sum_{n=1}^\infty K^n(t,s)K_n(t,s).
\]
Assume that ${\mathbf K}(t,t)$ converges for all
$t\in{\mathbb R}$. For every choice of $t_0\in{\mathbb R}$, the
map which to $f\in{\mathcal H}(K)$ associates the function
\[
\varphi(f)(t)=\sum_{n=1}^\infty K_n(t,t_0)f(t)^n
\]
maps ${\mathcal H}(K)$ into ${\mathcal H}({\mathbf K})$.
Furthermore,
\begin{equation}
\label{rue_des_francs_bourgeois}
\|\varphi(f)\|_{{\mathcal H}({\mathbf K})}^2
\le\sum_{n=1}^\infty
 K_n(t_0,t_0)\|f\|^{2n}_{{\mathcal H}(K)}.
\end{equation}
\end{Tm}

{\bf Proof:} Using \eqref{quartier_latin} and the two inequalities
in \eqref{saitoh_1} one after the other, we have
\[
\begin{split}
\|\varphi(f)\|_{{\mathcal H}({\mathbf K})}^2&\le\sum_{n=1}^\infty
\|K_n(\cdot, t_0)f(\cdot)^n\|^2_{{\mathcal H}(K_nK^n)}\\
&\le \sum_{n=1}^\infty\|K_n(\cdot, t_0)\|^2_{{\mathcal H}(K_n)}
\|f^n\|^2_{{\mathcal H}(K^n)}\\
&\le\sum_{n=1}^\infty
 K_n(t_0,t_0)\|f\|^{2n}_{{\mathcal H}(K)}
\end{split}
\]
since
\[
\|K_n(\cdot, t_0)\|^2_{{\mathcal H}(K_n)}=K_n(t_0,t_0).\]
\mbox{}\qed\mbox{}\\

As a consequence we have:

\begin{Tm} Let $t_0\in{\mathbb R}\setminus\{0\}$, and
let $\psi$ be the map which to $f\in{\mathcal H}(K_r)$ associates
the function
\begin{equation}
\label{lydia}
\psi(f)(t)=\sum_{n=1}^\infty\frac{\Gamma(n,r(t)+r(t_0))}{n!}f(t)^n.
\end{equation}
Then $\psi(f)$ belongs to the reproducing kernel Hilbert space
${\mathcal H}(K_{\psi,r})$ with reproducing kernel
$K_{\psi,r}(t,s)=\psi(r(t)+r(s))-\psi(r(t-s))$ and we have
\begin{equation}
\label{asie} \|\psi(f)\|_{{\mathcal
H}(K_{\psi,r})}^2\le\sum_{n=1}^\infty
\frac{\Gamma(n,2r(t_0))}{n!}\|f\|_{{\mathcal H}(K_r)}^{2n}
\end{equation}
\end{Tm}
Indeed, the kernel $\frac{\Gamma(n, r(t)+r(s))}{n!}$ 
is positive on the real line, and
\[
\left\|\frac{\Gamma(n, r(\cdot)+r(t_0))}{n!}\right\|^2_{{\mathcal H}\left(
\frac{\Gamma(n, r(t)+r(s))}{n!}\right)}=
\frac{\Gamma(n,2r(t_0))}{n!};
\]
See formula \eqref{poisse} if need be.\\

Similar theorems can be stated for the decompositions \eqref{St
Philippe du Roule, ligne 9} and \eqref{ohohoh}. The case
\eqref{ohohoh} will not be pursued here. To present the results
pertaining to the case \eqref{St Philippe du Roule, ligne 9}, we
make a number of additional assumptions. First we restrict
ourselves to positive $t$ and $s$. The kernel $(|t|+|s|-|t-s|)^n$
then takes the form
\[
K_{1,1}^n(t,s)=(|t|+|s|-|t-s|)^n=2^n\max(t^n, s^n),\]
and the corresponding reproducing kernel Hilbert space is equal to the
set of functions of the form
\[
F(t)=\int_0^{t^n} f(u)du, \quad f\in{\mathbf L}_2({\mathbb R}),\]
with norm (recall \eqref{poisse})
\[
\|F\|^2_{{\mathcal H}(K_{1,1}^n)}=
\frac{
\int_0^\infty |f|^2(u)du}{2^n}=\frac{\int_0^\infty
\frac{|F^\prime(u)|^2}{nu^{n-1}}du}{2^n}.\]
See \cite{MR1423330}, \cite{saitoh2}.
We will also assume that for some $t_0\ge 0$
\begin{equation}
\label{edf_gdf}
\sum_{n=1}^\infty \frac{\int_0^\infty u^ne^{-2ut_0} dm(u)}{2^n}<\infty.
\end{equation}

We have:

\begin{Tm} Assume that \eqref{edf_gdf} holds.
The map which to $F\in{\mathcal H}(K_{1,1})$ associates
\[
\psi(F)(t)=\sum_{n=1}^\infty\left(\int_0^{t^n}F^\prime(u)du\right)
\left (\int_0^\infty u^ne^{-u(t+t_0)}dm(u)\right)
\]
sends ${\mathcal H}(K_{1,1})$ into ${\mathcal H}(K_r)$ and we have:
\[
\|\psi(F)\|^2_{{\mathcal H}(K_r)}\le\|F^\prime\|^2_{{\mathbf
    L}_2({\mathbb R}^+)}\left\{
\sum_{n=1}^\infty \frac{\int_0^\infty u^ne^{-2ut_0} dm(u)}{2^n}\right\}.\]
\end{Tm}

In the case of the fractional Brownian motion, the above theorem
becomes (recall that $K_{H,1}(t,s)$ is defined by \eqref{ioanna}):

\begin{Tm} Let $0<H<1/2$ and let $t_0>1/2$.
The map which to $F\in{\mathcal H}(K_{1,1})$ associates the
\[
\psi(F)(t)=\frac{2H}{\Gamma(1-2H)}
\sum_{n=2}^\infty\left(\int_0^{t^n}F^\prime(u)du\right)
\frac{\Gamma(n-2H)}{n!(t+t_0)^{n-2H}}
\]
sends ${\mathcal H}(K_{1,1})$ into ${\mathcal H}(K_{H,1})$ and we have:
\begin{equation}
\label{totototo}
\|\psi(F)\|^2_{{\mathcal
H}(K_r)}\le\|F^\prime\|^2_{{\mathbf
    L}_2({\mathbb R}^+)}
\sum_{n=2}^\infty \frac{\Gamma(n-2H)}{2^nn!(2t_0)^{n-2H}}.
\end{equation}
\end{Tm}

We note that the series \eqref{totototo} converges  for $t_0>1/2$
since $\Gamma(n)=(n-1)!$.

\end{document}